%% file: robles_arXiv2.tex
\documentclass[12pt]{amsart}
\usepackage{amssymb,latexsym,amsmath,comment}
\usepackage[mathscr]{eucal}
\usepackage{bbm}
\usepackage{xcolor}

\numberwithin{equation}{section}
\numberwithin{table}{section}
\numberwithin{figure}{section}

\input{macros.tex}
\input{thms.tex}

\def\tNC{\mathcal{C}}

\def\bi{\hbox{\small{$\mathbf{i}$}}}

\def\CKpmhs{MR664326}

\def\CKSdeg{MR840721}

\def\CoMc{MR1251060}
\def\DeligneII{MR0498551}

\def\GGKtcu{MR3115136}

\def\Schmid{MR0382272}

\begin{document}
\title{Nilpotent cones and adjoint orbits}
\author[Robles]{C. Robles}
\email{robles@math.tamu.edu}
\address{Mathematics Department, Mail-stop 3368, Texas A\&M University, College Station, TX  77843-3368} 
\thanks{This work is partially supported by NSF grant DMS-1309238, and was undertaken while I was a member of the Institute for Advanced Study; I thank the institute for a wonderful working environment and the Robert and Luisa Fernholz Foundation for financial support.}
\date{\today}
\begin{abstract}
We show that the elements of the (open) cone underlying a nilpotent orbit on a period domain are pairwise congruent under (connected component of) the symmetry group of the period domain.
\end{abstract}
\keywords{Period domain, nilpotent orbit, \kahler~cone}
\subjclass[2010]
{
 14D07, 32G20. 
 58A14, 
}
\maketitle
\setcounter{tocdepth}{2}
\let\oldtocsection=\tocsection
\let\oldtocsubsection=\tocsubsection
\let\oldtocsubsubsection=\tocsubsubsection
\renewcommand{\tocsection}[2]{\hspace{0em}\oldtocsection{#1}{#2}}
\renewcommand{\tocsubsection}[2]{\hspace{3em}\oldtocsubsection{#1}{#2}}
\tableofcontents

\section{Introduction}

\subsection{Main result}

The purpose of this note is to give a simple proof of\footnote{See Section \ref{S:nt} for notation and definitions.}

\begin{theorem} \label{T:1}
Let $D$ be a period domain with symmetry group $G_\bR = \tAut(V_\bR,Q)$, and let $G_\bR^o \subset G_\bR$ be the connected component containing the identity.  Given a nilpotent orbit $(F^\sb;N_1,\ldots,N_m)$ on $D$, the underlying (open) nilpotent cone $\tNC\subset \fg_\bR$ is contained in an $\tAd(G_\bR^o)$--orbit.
\end{theorem}

Greg Pearlstein has observed that Theorem \ref{T:1} may be deduced from the several variables $\tSL_2$--orbit Theorem \cite[Theorem 4.20.iv]{\CKSdeg}.  However, to invoke so substantial and powerful a result could be misleading in this case: as I will demonstrate, Theorem \ref{T:1} follows from an elementary and geometrically transparent argument.

Nilpotent cones appear in a variety of contexts, and it is a longstanding problem to understand their structure.  Examples include: (i) the cone of \kahler~ classes on a smooth compact \kahler~ manifold, and (ii) the strictly convex conewise linear maps in the combinatorial intersection cohomology of polytopes; \cf~\cite{MR2429243} and the references therein.

An algebro--geometric variant of Theorem \ref{T:1} is addressed in \cite{GGVan}.  There one considers:
\begin{a_list}
\item
the stratification of the tangent space $T_X\mathrm{Def}(X)$ to the Kuranishi space when $X$ is a normal crossing divisor: very roughly, the strata correspond to the amount of smoothing that $X$ undergoes when deformed in the direction of the strata;
\item
the stratification of abelian subspaces of $\tNilp(\fg_\bR)$ by $\tAd(G_\bR)$--orbits.
\end{a_list}
The analog of Theorem \ref{T:1} is to show that the strata in (a) are mapped to strata in (b).

Theorem \ref{T:1} holds in the more general setting of Mumford--Tate domains \cite{SL2}.  To prove the theorem in this more general setting requires a fair amount of representation theory.  However, the argument in the case of period domains is comparatively simple, and captures the essential geometric idea:  the $\tAd(G_\bR)$--orbit of $N \in \tNC$ is characterized by data encoded in the limiting mixed Hodge structure $(F^\sb,W_\sb(N))$ associated with the nilpotent orbit $(F^\sb,N)$.  Theorem \ref{T:1} will follow from: (i) Springer and Steinberg's \cite{MR0268192} classification of the $\tAd(G_\bR^o)$--orbits in $\tNilp(\fg_\bR)$, and a result of Gerstenhaber and Hesselink \cite{MR0429875, MR0136683} on the partial order on these orbits; and (ii) Cattani and Kaplan's \cite{\CKpmhs} result on the independence of the weight filtration $W_\sb(N)$ on our choice of $N \in \tNC$.

The results of Springer and Steinberg, and Gerstenhaber and Hesselink apply to the case that $G_\bR$ is a classical Lie group.  To prove the analog of Theorem \ref{T:1} in the more general setting of a Mumford--Tate domain $\sD = \sG_\bR/R$ requires:
\begin{i_list}
\item The \dokovic--Kostant--Sekiguchi \cite{MR891636, MR867991} classification of the nilpotent $\tAd(\sG_\bR)$--conjugacy classes for an arbitrary (connected) semisimple real Lie group $\sG_\bR$.
\item
The observation that underlying any nilpotent orbit is a Hodge--Tate degeneration \cite{GGR}, and the classification of the latter.
\end{i_list}
The components of (ii) yield classifications of the horizontal $\tSL_2$'s and the $\bR$--split limiting mixed Hodge structures on a Mumford--Tate domain domain, up to the action of $\sG_\bR$, \cite{SL2}.  This is related to work in progress, by Matt Kerr, Greg Pearlstein and the author, to understand the relationship between 
\begin{bcirclist}
\item 
  $\sG_\bR$--orbits of limiting mixed Hodge structures, and 
\item
  $\sG_\bR$--orbits in the closure $\tcl(\sD)$ of $\sD$ in the compact dual $\sG_\bC/P$.
\end{bcirclist}  
These two sets of $\sG_\bR$--orbits are related by the $\sG_\bR$--equivariant reduced/\naive~ limit period mapping $\Phi_\infty : B(N) \to \tcl(\sD)$, \cite{KP2013, \GGKtcu, GGR}; here $B(N)$ is the boundary component associated with the nilpotent $N \in \fg_\bR$.  The image $\Phi_\infty(B(N))$ lies in a $\sG_\bR$--orbit $\sO \subset \tcl(\sD)$.  More generally, the image under $\Phi_\infty$ of the $\sG_\bR$--orbit of a $\bR$--split limiting mixed Hodge structure $(F^\sb,N)$ is equal to a $\sG_\bR$--orbit $\sO \subset \tcl(\sD)$.  This yields a map
\[
  \{\hbox{$\sG_\bR$--orbits of $\bR$--split LMHS}\} \ \to \ 
  \{\hbox{$\sG_\bR$--orbits in $\tcl(\sD)$}\} \,;
\]
one goal is to determine the image and fibres of the map, and the effect of the map on the natural partial orders (given by containment in closure) on each set.

We finish the introduction by reviewing some notation and results from Hodge theory.  The necessary material from representation theory is reviewed in Section \ref{S:repthy}, and Theorem \ref{T:1} is proved in Section \ref{S:prf}.

\subsection*{Acknowledgements}
Over the course of this work I have benefited from conversations and correspondence with a number of colleagues.  I would especially like to thank Eduardo Cattani, Mark Green, Phillip Griffiths, Matt Kerr and Greg Pearlstein for their time and insight.

\subsection{Notation and terminology} \label{S:nt}

Fix a real vector space $V_\bR$.  A (pure, real) \emph{Hodge structure of weight $k$} on $V_\bR$ is a decomposition
\[
  V_\bC \ = \ \bigoplus_{p+q=k} V^{p,q}
\]
of the complexification $V_\bC = V_\bR \ot_\bR \bC$ with the property
\[
  \overline{V^{p,q}} \ = \ V^{q,p} \,.
\]
The \emph{Hodge numbers} are $\bh = (h^{p,q} = \tdim_\bC\,V^{p,q})$.  The \emph{Hodge flag} $F^\sb$ associated with the Hodge structure is the filtration of $V_\bC$ given by
\[
  F^p \ = \ \bigoplus_{r\ge p} V^{r,\sb} \,.
\]
We may recover the Hodge structure from the Hodge flag by 
\[
  V^{p,q} \ = \ F^p \,\cap\, \overline{F^{q}} \,.
\]

Fix a nondegenerate $(-1)^k$--symmetric bilinear form $Q$ on $V_\bR$.  The Hodge structure is \emph{$Q$--polarized} if
\[
  Q(V^{p,q},V^{r,s}) \ = \ 0 \quad\hbox{when}\quad (p,q) \not= (s,r)
\]
and 
\[
  \bi^{p-q} Q(v,\bar v) \ > \ 0 
  \quad\hbox{for all}\quad 0\not= v \in V^{p,q}\,,
\]
where $\bi = \sqrt{-1}$.  The \emph{period domain} $D = D_\bh$ is the set of all $Q$--polarized Hodge structures on $V_\bR$ with Hodge numbers $\bh$; is a $G_\bR = \tAut(V_\bR,Q)$--homogeneous complex manifold.  

The Hodge flag is $Q$--isotropic
\[
  Q(F^p,F^{n+1-p}) \ = \ 0 \,,
\]
and thus an element of the \emph{compact dual}, the $G_\bC = \tAut(V_\bC,Q)$--homogeneous flag variety
\[
  \check D \ = \ \tFlag^Q_\mathbf{f} V_\bC  
\]
consisting of $Q$--isotropic filtrations $F^\sb$ of $V_\bC$ with dimension $\mathbf{f} = (f^p = \tdim_\bC F^p)$.  Taking the $G_\bR$--orbit of the Hodge flag $F^\sb\in\check D$ realizes the period domain $D$ as an open submanifold of the compact dual.

Let $\fg_\bR = \tEnd(V_\bR,Q)$ be the Lie algebra of $G_\bR$.   A ($m$--variable) \emph{nilpotent orbit} on $D$ consists of a tuple $(F^\sb; N_1,\ldots,N_m )$ such that $F^\sb \in \check D$, the $N_i \in \fg_\bR$ commute and $N_iF^p \subset F^{p-1}$, and the holomorphic map $\psi : \bC^m \to \check D$ defined by
\begin{equation}\label{E:htno}
  \psi(z^1,\ldots,z^m) \ = \ \texp( z^i N_i ) F^\sb 
\end{equation}
has the property that $\psi(z) \in D$ for $\tIm(z^i) \gg 0$.  The associated (open) \emph{nilpotent cone} is
\begin{equation} \label{E:s}
  \tNC\ = \ \{ t^i N_i \ | \ t^i > 0 \} \,.
\end{equation}

Given a nilpotent $N \in \fg_\bR$ such that $N^{k+1}=0$, there exists a unique increasing filtration $W_0(N) \subset W_1(N) \subset \cdots \subset W_{2k}(N)$ of $V_\bR$ with the properties that 
\[
  N\,W_\ell(N) \ \subset \ W_{\ell-2}(N)
\]
and the induced 
\[
  N^\ell : \tGr_{k+\ell} W_\sb(N) \ \to \ \tGr_{k-\ell} W_\sb(N)
\]
is an isomorphism for all $\ell \le k$.  Above, $\tGr_m W_\sb(N) = W_m(N)/W_{m-1}(N)$.  Moreover, 
\[
  Q_\ell(u,v) \ = \  Q(u , N^\ell v)
\]
defines a nondegenerate $(-1)^{k+\ell}$--symmetric bilinear form on $\tGr_{k+\ell} W_\sb(N)$.  

Define
\[
  \tGr_{k+\ell} W_\sb(N)_\tprim \ = \ 
  \tker\,\{
    N^{\ell+1} : \tGr_{k+\ell} W_\sb(N) \to \tGr_{k-\ell-2} W_\sb(N) \} \,,
\]
for all $\ell \ge 0$.  A \emph{limiting mixed Hodge structure} (or \emph{polarized mixed Hodge structure}) on $D$ is given by a pair $(F^\sb,N)$ such that $F^\sb \in \check D$, $N \in \fg_\bR$ and $N(F^p) \subset F^{p-1}$, the filtration $F^\sb$ induces a weight $m$ Hodge structure on $\tGr_m W_\sb(N)$ for all $m$, and the Hodge structure on $\tGr_{k+\ell}(W_\sb(N))_\tprim$ is polarized by $Q_\ell$ for all $\ell\ge0$.  The notions of nilpotent orbit and limiting mixed Hodge structure closely related.  Indeed, they are equivalent when $m=1$.

\begin{theorem}[Cattani, Kaplan, Schmid] \label{T:cks}
Let $D \subset \check D$ be a Mumford--Tate domain (and compact dual) for a Hodge representation of $G_\bR$.
\begin{a_list_emph}
\item 
A pair $(F^\sb;N)$ forms a one--variable nilpotent orbit if and only if it forms a limiting mixed Hodge structure, \emph{\cite[Corollary 3.13]{\CKSdeg} and \cite[Theorem 6.16]{\Schmid}}.  
\item 
Given an $m$--variable nilpotent orbit $(F^\sb;N_1,\ldots,N_m)$, the weight filtration $W_\sb(N)$ does not depend on the choice of $N\in\tNC$, \emph{\cite[Theorem 3.3]{\CKpmhs}.  Let $W_\sb(\tNC)$ denote this common weight filtration.} 
\end{a_list_emph}
\end{theorem}

The \emph{Deligne splitting} \cite{\CKSdeg, \DeligneII} 
\begin{subequations} \label{SE:deligne}
\begin{equation}
  V_\bC \ = \ \bigoplus I^{p,q}
\end{equation}
associated with the limiting mixed Hodge structure is given by
\begin{equation}
  I^{p,q} \ = \ F^p \,\cap\, W_{p+q} \, \cap \, 
  \Big( \overline{F^q} \,\cap\,W_{p+q} \,+\, 
         \sum_{j\ge1} \overline{F^{q-j}} \,\cap\, W_{p+q-j-1} \Big) \,.
\end{equation}
\end{subequations}
It is the unique bigrading of $V_\bC$ with the properties that 
\begin{equation} \label{E:FW}
  F^p \ = \ \bigoplus_{r \ge p} I^{r,\sb} \tand
  W_\ell(\tNC) \ = \ \bigoplus_{p+q \le \ell} I^{p,q} \,,
\end{equation}
and
\[
  \overline{I^{p,q}} \ = \ I^{q,p} \quad\hbox{mod} \quad
  \bigoplus_{r<q,s<p} I^{r,s} \,.
\]
Observe that $N I^{p,q} \subset I^{p-1,q-1}$ for all $N \in \tNC$.  Given $\ell \ge0$ and $p+q = k+\ell$, we define
\[
  I^{p,q}_\tprim \ = \ \tker\,\{
  N^{\ell+1} : I^{p,q} \to I^{-p-1,-q-1} \} \,.
\]

\section{Representation theory background} \label{S:repthy}

This is a terse summary of the classification of nilpotent conjugacy classes in symplectic and orthogonal Lie algebras.  Collingwood and McGovern's \cite{\CoMc} is an excellent reference for the material in this section.\footnote{The term `nilpotent orbit' is used in both Hodge theory and representation theory to refer to two distinct (but related) objects.  In \cite{\CoMc} a `nilpotent orbit' is the $\tAd(G)$--orbit of an nilpotent element of the Lie algebra $\fg$.  To avoid confusion I will refer to the representation theoretic version as a `nilpotent conjugacy class.'}

Let $\bbk$ be one of $\bR,\bC$.  Let $V_\bbk$ be a finite--dimensional $\bbk$--vector space.  A \emph{standard triple} is a set $\{ M , Y , N \} \subset \tEnd(V_\bbk)$ with the property
\begin{equation} \label{E:st}
  [Y , M] = 2M \,,\quad [Y,N] = -2N \tand [M,N] = Y \,.
\end{equation}
We say that $N$ is the \emph{nilnegative element} of the standard triple.  From \eqref{E:st} we see that $\tspan\{M,Y,N\}\simeq\fsl_2\bbk$ is a three--dimensional semisimple subalgebra (TDS) of $\tEnd(V_\bbk)$.

Fix a nilpotent element $N \in \tEnd(V_\bbk)$.  The Jacobson--Morosov Theorem \cite{MR0049882, MR0007751} asserts that $N$ is the nilnegative element of a standard triple.  Let $\fsl_2\bbk \subset \tEnd(V_\bbk)$ be the TDS spanned by a standard triple containing $N$ as the nilnegative element, and let
\begin{equation} \label{E:V}
  V_\bbk \ = \ \bigoplus_{\ell\ge0} V(\ell)
\end{equation}
be the $\fsl_2\bbk$--decomposition of $V_\bbk$.  Here $V(\ell)$ is the \emph{isotypic component of highest weight $\ell$}; it is isomorphic to the direct sum $(\tSym^\ell\bbk^2)^{\op m_\ell}$ of $m_\ell$ irreducible $\fsl_2\bbk$--modules of dimension $\ell+1$.  In particular, $V(\ell)$ admits a basis of the form 
\[
  \{ N^a v_i \ | \ 1 \le i \le m_\ell \,,\ 0 \le a \le \ell \} \,.
\]
Here $N^\ell v_i \not=0$ and $N^{\ell+1}v_i = 0$.  Loosely, we think of $V(\ell)$ as spanned by $m_\ell$ $N$--strings of length $\ell+1$.   The Jordan normal form for elements of $\tEnd(V_\bbk)$ implies that two nilpotents $N , N' \in \tEnd(V_\bbk)$ lie in the same $\tAut(V_\bbk)$--orbit (under the adjoint action) if and only if $m_\ell = m_\ell'$ for all $0 \le \ell \in \bZ$.

The \emph{highest weight space of $V(\ell)$} is  
\begin{equation} \label{E:P(l)}
  P(\ell) \ = \ \tspan_\bR \{ v_i \ | \ 1 \le i \le m_\ell\} \,.
\end{equation}
Note that
\[
  V(\ell) \ = \ \bigoplus_{a=0}^\ell N^a P(\ell) \,.
\]
The weight filtration of Section \ref{S:nt} is 
\[
  W_\ell(N) \ = \ \bigoplus_{m-2a \le \ell-k} N^a P(m) \,.
\]
Observe that 
\begin{equation} \label{E:PvI}
  P(\ell) \ \simeq \ \tGr_{\ell+k} W_\sb(N)_\tprim 
  \simeq \ \bigoplus_{p+q=\ell+k} I^{p,q}_\tprim \,.
\end{equation}

Let $Q$ be a nondegenerate $(-1)^k$--symmetric bilinear form on $V_\bbk$ and let $G_\bbk = \tAut(V_\bbk,Q)$ be the corresponding Lie group with Lie algebra $\fg_\bbk = \tEnd(V_\bbk,Q)$.  Given a nonzero $N \in \tNilp(\fg_\bbk)$, the Jacobson--Morosov Theorem asserts that $N$ is the nilnegative element of a standard triple contained in $\fg_\bbk$.  So, we may assume that the TDS is contained in $\fg_\bbk$.  Then 
\[
  Q_\ell(u,v) \ = \ Q( u , N^\ell v )
\]
defines a non--degenerate $(-1)^{k+\ell}$--symmetric bilinear form on $P(\ell)$.  (This implies $m_\ell$ is even if $k+\ell$ is odd.)  

Let $\tNilp(\fg_\bbk) \subset \fg_\bbk$ denote the set of nilpotent elements.  Then $\tNilp(\fg_\bbk)$ decomposes into a finite union of $\tAd(G_\bbk)$--orbits.  Each such orbit is a \emph{nilpotent conjugacy class}.  It is clear that $\bm = (m_0,\ldots,m_k)$ is an invariant of the nilpotent conjugacy class $\cN$ of $N$ which we will denote $\bm_\cN$.  In the case that $V$ and $Q$ are defined over $\bC$, $\bm_\cN$ is a complete invariant of the nilpotent conjugacy class.  More precisely, we have:

\begin{theorem}[{Gerstenhaber \cite{MR0136683}}] \label{T:C}
Let $Q$ be a $(-1)^k$--symmetric bilinear form on a complex vector space $V_\bC$, and set $G_\bC = \tAut(V_\bC,Q)$ with Lie algebra $\fg_\bC = \tEnd(V_\bC,Q)$.  The $\tAd(G_\bC)$--conjugacy classes $\cN\subset \tNilp(\fg_\bC)$ are uniquely determined by the $\bm_\cN$.
\end{theorem}

\noindent Theorem \ref{T:C} does not quite hold when we replace $G_\bC = \tAut(V_\bC,Q)$ with the connected component $G_\bC^o \subset G_\bC$ containing the identity.  The difficulty is that a single $\tAd(G_\bC)$--orbit $\cN \subset \tNilp(\fg_\bC)$ may decompose into two $\tAd(G_\bC^o)$--orbits.

\begin{theorem}[{Gerstenhaber \cite{MR0136683}, Springer--Sternberg \cite{MR0268192}}] \label{T:Co}
Let $Q$ be a $(-1)^k$--symmetric bilinear form on a complex vector space $V_\bC$, and let $G_\bC^o$ be the connected component of $\tAut(V_\bC,Q)$ containing the identity.  Suppose that $\bm = \bm_\cN$ for some $\tAd(G_\bC)$--conjugacy class $\cN \subset \tNilp(\fg_\bC)$.
\begin{a_list_emph}
\item
Suppose that $Q$ is symmetric, $m_{2\ell}=0$ and $m_{2\ell+1}$ is even for all $\ell$.  Then there exist exactly two distinct nilpotent $\tAd(G_\bC^o)$--conjugacy classes with this invariant ($\cN$ decomposes into two $\tAd(G^o_\bC)$--orbits).
\item 
Otherwise, $\cN$ is the unique $\tAd(G_\bC^o)$--orbit with this invariant.
\end{a_list_emph}
\end{theorem}

Suppose that $\bbk = \bR$.  Let $s_\ell = (p_\ell,q_\ell)$ be the signature of $Q_\ell$, when $k+\ell$ is even.  Set $\bs = (s_\ell)$.  It is clear that $\bs$ is an invariant of the nilpotent conjugacy class $\cN \subset \tNilp(\fg_\bR)$ which we will denote by $\bs_\cN$.  The $\tAd(G_\bR)$--conjugacy classes $\cN \subset \tNilp(\fg_\bR)$ are classified by the $\bm$ and $\bs$.  More precisely, we have:  

\begin{theorem}[{Burgoyne--Cushman \cite{MR0432778}, Springer--Steinberg \cite{MR0268192}}] \label{T:R}
Let $Q$ be a nondegenerate $(-1)^k$--symmetric bilinear form on a real vector space $V_\bR$, and set $G_\bR = \tAut(V_\bR,Q)$ with Lie algebra $\fg_\bR = \tEnd(V_\bR,Q)$.  The $\tAd(G_\bR)$--conjugacy classes $\cN\subset \tNilp(\fg_\bR)$ are uniquely determined by the $\bm_\cN$ and $\bs_\cN$.
\end{theorem}

\noindent  Let $G_\bR^o \subset \tAut(V_\bR,Q)$ denote the connected component containing the identity.  As above, the $\tAd(G_\bR^o)$--conjugacy classes $\cN \subset \tNilp(\fg_\bR)$ are \emph{almost} classified by the $\bm$ and $\bs$.    

\begin{theorem}[{Burgoyne--Cushman \cite{MR0432778}, Springer--Steinberg \cite{MR0268192}}] \label{T:Ro}
Let $Q$ be a nondegenerate $(-1)^k$--symmetric bilinear form on a real vector space $V_\bR$, and let $G_\bR^o \subset \tAut(V_\bR,Q)$ be the connected component containing the identity with Lie algebra $\fg_\bR = \tEnd(V_\bR,Q)$.  Suppose that $\bm = \bm_\cN$ and $\bs=\bs_\cN$ are invariants for some $\tAd(G_\bR)$--conjugacy class $\cN \subset \tNilp(\fg_\bR)$.  Suppose that $Q$ is symmetric and that one of the following conditions holds:
\begin{a_list}
\item $m_{2\ell}=0$ for all $\ell$ (so the nontrivial $Q_\ell$ are all skew--symmetric), 
\item some $m_{2\ell}\not=0$, and either: \emph{(i)} $Q_{2\ell}$ is $(-1)^\ell$--definite for every $m_{2\ell} \not=0$, or \emph{(ii)} $Q_{2\ell}$ is $(-1)^{\ell+1}$--definite for every $m_{2\ell} \not=0$.
\end{a_list}
Then $\cN$ decomposes into two $\tAd(G_\bR^o)$--orbits.  Otherwise $\cN$ is an $\tAd(G_\bR^o)$--orbit.
\end{theorem}

\begin{remark}[Orbit closure] \label{R:orbcl}
Let $\cN$ be an $\tAd(G_\bR^o)$--orbit in the cone $\tNilp(\fg_\bR)$ of nilpotent elements.  Then the closure of $\cN$ is a union 
\[
  \overline{\cN} \ = \ \bigcup_{\cN' \subset \overline\cN} \cN' 
\]
$\tAd(G_\bR^o)$--orbits $\cN' \subset \tNilp(\fg_\bR)$.  Define a partial order on the tuples $\bm$ by declaring $\bm \ge \bm'$ if 
\[
  \sum_{1 \le i \le \ell} m_i \ \ge \ \sum_{1 \le i \le \ell} m_i'
  \quad\hbox{for all} \quad \ell \,. 
\]
If $\cN' \subsetneq \overline\cN$ then it is necessarily the case that $\bm' < \bm$, \cite{MR0429875, MR0136683}.
\end{remark}

\section{Proof of Theorem \ref{T:1}} \label{S:prf}

Theorem \ref{T:1} is a consequence of Theorems \ref{T:cks}(b) and \ref{T:Ro}, and Remark \ref{R:orbcl}.  Suppose that $N,N' \in \tNC$.  Theorem \ref{T:cks}(b) implies that $W_\sb(N) = W_\sb(N')$.  Then \eqref{E:PvI} and the polarization conditions imply  
\begin{equation}\label{E:1}
  \bm(N) \ = \ \bm(N') 
  \tand 
  \bs(N) \ = \ \bs(N') \,.
\end{equation}
If we are in case (c) of Theorem \ref{T:R}, then Theorem \ref{T:1} follows directly.

If we are not in case (c), then it is possible that the $\tAd(G_\bR^o)$--conjugacy classes $\cN$ and $\cN'$ of $N$ and $N$ are distinct.  In this case Remark \ref{R:orbcl} and \eqref{E:1} imply 
\[
  \cN \,\cap\, \overline\cN{}' \ = \ \emptyset \ = \ 
  \cN' \,\cap\, \overline\cN \,.
\]
This, along with the fact that the nilpotent cone $\tNC$ is open and connected, forces $\cN = \cN'$.
\hfill\qed

\def\cprime{$'$} \def\Dbar{\leavevmode\lower.6ex\hbox to 0pt{\hskip-.23ex
  \accent"16\hss}D}

\end{document}

%% file: macros.tex
\def\kahler{{K\"ahler}}
\def\naive{{na\"ive}}

\def\dokovic{{\DH}okovi{\'c}}
\def\cf{cf.~}

\newcommand{\hsp}[1]{{\hbox{\hspace{#1}}}}

\def\tAd{\mathrm{Ad}} 
 
\def\tAut{\mathrm{Aut}}

\def\bC{\mathbb C}

 \def\tcl{\mathrm{cl}}
 
 \def\sD{\mathscr{D}}

 \def\tdim{\mathrm{dim}}

\def\tEnd{\mathrm{End}}

\def\texp{\mathrm{exp}}

\def\tFlag{\mathrm{Flag}}

\def\sG{\mathscr{G}}

\def\tGr{\mathrm{Gr}}
\def\fg{{\mathfrak{g}}}

\def\bh{\mathbf{h}}
\def\bi{\mathbf{i}}

 \def\tIm{\mathrm{Im}}

\def\bbk{\mathbbm{k}}
 
 \def\tker{\mathrm{ker}}

 \def\bm{\mathbf{m}}

 \def\cN{\mathcal N}

\def\tNilp{\mathrm{Nilp}}

\def\sO{\mathscr{O}}

\def\tprim{\mathrm{prim}}

\def\bR{\mathbb R}

\def\bs{\mathbf{s}} 

\def\tSL{\mathrm{SL}}

 \def\tSym{\mathrm{Sym}}

 \def\tspan{\mathrm{span}}
\def\fsl{\mathfrak{sl}}

   \def\bZ{\mathbb Z}

\def\tand{\quad\hbox{and}\quad}

\def\sb{{\hbox{\tiny{$\bullet$}}}}

\def\op{\oplus}
\def\ot{\otimes}

\newcounter{numcnt}

\newcounter{cnt}

\newcounter{acnt}
\newenvironment{a_list}{ 
  \begin{list}{{(\alph{acnt})}}
   {\usecounter{acnt} \setlength{\itemsep}{3pt}
    \setlength{\leftmargin}{25pt} \setlength{\labelwidth}{20pt} }
   }
   {\end{list}}
\newenvironment{a_list_emph}{ 
  \begin{list}{{\emph{(\alph{acnt})}}}
   {\usecounter{acnt} \setlength{\itemsep}{3pt}
    \setlength{\leftmargin}{25pt} \setlength{\labelwidth}{20pt} }
   }
   {\end{list}}

\newcounter{Acnt}

\newcounter{icnt}

\newenvironment{i_list}{ 
  \begin{list}{{(\roman{icnt})}}
   {\usecounter{icnt} \setlength{\itemsep}{3pt}
    \setlength{\leftmargin}{25pt} \setlength{\labelwidth}{20pt} }
   }
   {\end{list}}
\newcounter{Icnt}

\newcounter{exam_cnt}

\newcounter{mccnt}

\newenvironment{bcirclist}{ 
  \begin{list}{\boldmath$\circ$\unboldmath}
   {\usecounter{cnt} \setlength{\itemsep}{2pt}
    \setlength{\leftmargin}{15pt} \setlength{\labelwidth}{20pt} }
   }
   {\end{list}}

\def\CoMc{MR1251060}


%% file: thms.tex
\newtheorem*{corollary*}{Corollary}

\newtheorem*{lemma*}{Lemma}

\newtheorem*{proposition*}{Proposition}
\newtheorem{theorem}[equation]{Theorem}
\newtheorem*{theorem*}{Theorem}

\theoremstyle{definition}

\newtheorem*{boldQ*}{Question}
\newtheorem*{boldP*}{Problem}

\theoremstyle{remark}
\newtheorem*{assume*}{Assume}
\newtheorem*{answer*}{Answer}

\newtheorem*{claim*}{Claim}

\newtheorem*{definition*}{Definition}

\newtheorem*{example*}{Example}
\newtheorem*{hint*}{Hint}
\newtheorem*{notation*}{Notation}
\newtheorem{remark}[equation]{Remark}
\newtheorem*{remark*}{Remark}
\newtheorem*{remarks*}{Remarks}
\newtheorem*{fact*}{Fact}
\newtheorem*{emphL*}{Lemma}

\newtheorem*{emphQ*}{Question}
\newtheorem*{emphA*}{Answer}

\numberwithin{HWeq}{section}
\theoremstyle{definition}

